\documentclass[a4paper,11pt]{amsart}

\usepackage{amssymb,amscd,amsthm,mathrsfs}
\usepackage[ansinew]{inputenc}
\usepackage[centertags]{amsmath}

   \def\DD{{\mathbb D}}

 \def\NN{{\mathbb N}}  
\def\QQ{{\mathbb Q}} \def\RR{{\mathbb R}}  \def\TT{{\mathbb T}}
   
 \def\ZZ{{\mathbb Z}}

\def\cD{\mathcal{D}}    
    \def\cW{\mathcal{W}}
\def\cF{\mathcal{F}}

\newtheorem*{teo*}{Theorem}

\newtheorem{teo}{Theorem}[section]

\newtheorem{quest}{Question}

\newtheorem*{af}{Claim}
\newtheorem{lema}[teo]{Lemma}
\newtheorem{prop}[teo]{Proposition}

\newcommand{\bi}{\begin{itemize}}
\newcommand{\ei}{\end{itemize}}

\theoremstyle{definition}
\newtheorem*{defi}{Definition}
\theoremstyle{remark}
\newtheorem{obs}[teo]{Remark}

\newcommand{\dem}{\vspace{.05in}{\sc\noindent Proof.}}
\newcommand{\lqqd}{\par\hfill {$\Box$} \vspace*{.05in}}

\newcommand{\eps}{\varepsilon}

\newcommand{\en}{\subset}

 \DeclareMathOperator{\Emb}{Emb}

\DeclareMathOperator{\diametro}{diam}
\DeclareMathOperator{\Image}{Image}

\author[R. Potrie]{Rafael Potrie}
\address{CMAT, Facultad de Ciencias, Universidad de la Rep\'ublica, Uruguay}
\urladdr{www.cmat.edu.uy/$\sim$rpotrie}
\email{rpotrie@cmat.edu.uy}

\title{Partially hyperbolic diffeomorphisms with a trapping property}
\thanks{The autor was partially supported by CSIC group 618, FCE-3-2011-6749 and the Palis Balzan project.}

\begin{document}
\maketitle

\begin{abstract}
We study partially hyperbolic diffeomorphisms satisfying a
trapping property which makes them look as if they were Anosov at
large scale. We show that, as expected, they share several
properties with Anosov diffeomorphisms. We construct an expansive
quotient of the dynamics and study some dynamical consequences
related to this quotient.

%\marginpar{This note will appear in the margin.}

%\bigskip

%\noindent
%{\bf Keywords:}

%\medskip

%\noindent {\bf MSC 2000:} 37C05, 37C20, 37C25, 37C29, 37D30.
\end{abstract}

\section{Introduction}\label{SectionIntroduccion}

The purpose of this paper is twofold. On the one hand, we provide
some mild contributions to the classification problem of partially
hyperbolic diffeomorphisms in higher dimensions. On the other hand
we study the dynamics of certain partially hyperbolic
diffeomorphisms as well as give evidence of certain pathological
phenomena which must be dealt with in order to understand better
this kind of dynamics.

We shall consider partially hyperbolic diffeomorphisms $f: M \to
M$ admitting a splitting of the form $TM = E^{cs} \oplus E^u$ in
the pointwise sense (see below for precise definitions).

It is well known (\cite{HPS}) that the unstable bundle $E^u$ is
uniquely integrable into a foliation called \emph{unstable
foliation} and denoted as $\cW^u$. On the other hand, the
integrability of the bundle $E^{cs}$  is a subtler issue (see
\cite{BurnsWilkinson}). As in other results aiming at the
classification of partially hyperbolic diffeomorphisms in higher
dimensions (for example \cite{Bonhet,Carrasco,Gogolev}) we shall
ignore this issue for the moment by assuming that $f$ is
\emph{dynamically coherent}. We say that a partially hyperbolic
diffeomorphism is \emph{dynamically coherent} if there exists an
$f$-invariant foliation $\cW^{cs}$ tangent to $E^{cs}$. In
dimension 3, there are several classification type results which
do not assume integrability to start with (see for example
\cite{BonattiWilkinson,BBI,HammerlindlPotrie-Nil,HammerlindlPotrie-Solv,Pot-PHFoliations}).
See also section \ref{Section-Coherence} of this paper where we
remove this hypothesis under a different assumption.

This work will concern partially hyperbolic diffeomorphisms which verify a dynamical condition which makes them look, from far apart, as Anosov diffeomorphisms. We will say that a dynamically  coherent partially hyperbolic diffeomorphism $f: M \to M$ with splitting $TM = E^{cs} \oplus E^u$ has \emph{trapping property}
if there exists a continuous map $\cD^{cs}: M \to \Emb^1(\overline{\DD}^{cs}, M)$ such that $\cD^{cs}(x)(0)=x$, the image of $\overline{\DD}^{cs}$ (the closed unit ball of dimension $\dim E^{cs}$) by $\cD^{cs}(x)$ is always contained in $\cW^{cs}(x)$ and they verify the following trapping property:
$$f(\cD^{cs}(x)(\overline{\DD}^{cs})) \en \cD^{cs}(f(x))(\DD^{cs})
\qquad \forall x\in M.$$

The main point of this paper is to recover the same type of results that are valid for Anosov diffeomorphisms (\cite{Franks,Manning,Newhouse}) in this setting. Several examples enjoy this trapping property and it is important in order to obtain dynamical consequences (\cite{ManheContributions,Carvalho,BV,Pot-WildMilnor,BF,FisherPotrieSambarino,Roldan}), however, the point here is to avoid the usual assumption that the trapping property is seen at ``small scale'' (a similar approach is pursued in \cite{Pot-JDDE} in the 3-dimensional case and with one dimensional center). See also \cite{CP} for a related notion of \emph{chain-hyperbolicity}.

A relevant point is to obtain results concerning partially hyperbolic diffeomorphisms without knowledge \emph{a priori} on the global structure of the invariant foliations. To support this point of view we give in section \ref{Section-Coherence} a weaker assumption which implies dynamical coherence and our trapping
property.

\subsection{Statement of results}
In this paper, we will consider partially hyperbolic diffeomorphisms in one of the weakest forms (see \cite[Appendix B]{BDV} for a survey of possible definitions). We make explicit the definition we shall use to avoid confusions with other references.

A $C^1$-diffeomorphism $f: M \to M$ is \emph{partially hyperbolic} if there exists a $Df$-invariant continuous splitting $TM = E^{cs}\oplus E^u$ and $N \geq 1$ such that for every $x\in M$ and for every pair of unit vectors $v^{cs} \in E^{cs}(x)$ and $v^{u} \in E^u(x)$ one has that:

$$  \|Df^N v^u \| > \max \{ 1, \|Df^N v^{cs} \| \} $$

Notice that with this definition $f$ might be partially hyperbolic and $f^{-1}$ not. This will not be a problem here, the results can be easily adapted to other settings. The definitions of \emph{dynamical coherence} and the \emph{trapping property} are the ones given in the introduction.

\begin{teo} Let $f: M \to M$ be a dynamically coherent partially hyperbolic diffeomorphism with splitting $TM = E^{cs} \oplus E^u$ and satisfying the trapping property. Assume moreover that one of the following holds:
\begin{itemize}
\item $M= \TT^d$ or, \item $\dim E^u=1$.
\end{itemize}
\noindent Then, $M= \TT^d$ and $f$ is homotopic to a linear Anosov
automorphism of $\TT^d$.
\end{teo}

An expansive quotient of the dynamics will be constructed under the hypothesis of the theorem. This will be enough to obtain that results of non-existence of Anosov diffeomorphisms depending on the Lefschetz formula hold for partially hyperbolic diffeomorphisms in our setting (see for example \cite{NoAnosovEnAlgunasVariedades,GogolevHertz}). This is done in section \ref{Section-ExpansiveQuotient} and the end of section \ref{Section-TransitiveExpansive}. It might be that this quotient is of independent interest.

Under certain assumptions (resembling those of the theory of Franks-Newhouse-Manning \cite{Franks,Newhouse,Manning}) we will see in section \ref{Section-TransitiveExpansive} that the quotient map is in fact transitive by translating the proofs in the Anosov setting to ours.

In section \ref{Section-Dynamics} some mild dynamical consequences are derived and some questions which might help to understand better the panorama are posed. Finally, in section \ref{Section-Coherence} we give weaker conditions which ensure dynamical coherence.

Appendix \ref{Appendix-Decomposition} presents an example which improves a decomposition constructed in \cite{RobertsDescomposiciones} showing that the quotients might be quite wild while respecting (some of) the dynamical conditions.

\smallskip
{\bf Acknowledgements: }\textit{This is an improvement of part of my thesis \cite[Section 5.4]{Pot-Tesis} under S. Crovisier and M. Sambarino, I thank both for several discussions and motivation. I also benefited from discussions with A. Gogolev, N. Gourmelon and M. Rold\'an as well as important input by the referee. This work is dedicated to Jorge Lewowicz (1937-2014) who in particular has inspired several of the ideas here presented. }

\section{Notation}

Along this paper $M$ will denote a closed $d$-dimensional manifold and $f: M \to M$ a partially hyperbolic diffeomorphism. Except in section \ref{Section-Coherence} we shall assume that $f$ is dynamically coherent and verifies the trapping property. Along the paper we shall assume that $d\geq 3$ (for the case of $d=2$ stronger results can be obtained with easier proofs, see for example \cite[Section 4.A]{Pot-Tesis}).

Given any foliation $\cF$ on $M$ we shall denote as $\cF(x)$ to the leaf through $x$, $\cF_{\eps}(x)$ to the $\eps$-disk around $x$ in the induced metric of the leaf and $\tilde \cF$ will always denote the lift of $\cF$ to the universal cover $\tilde M$ of $M$. Here, foliation means a continuous foliation with $C^1$-leafs tangent to a continuous distribution (foliations of class $C^{1,0+}$ according to \cite{CandelConlon}).

\section{An expansive quotient of the
dynamics}\label{Section-ExpansiveQuotient}

Denote as $\overline{\cD}^{cs}_x$ to
$\cD^{cs}(x)(\overline{\DD}^{cs})$ and $\cD^{cs}_x$ to
$\cD^{cs}(x)(\DD^{cs})$.

We can define for each $x\in M$

$$ A_x = \bigcap_{n\geq 0} f^{n}(\overline{\cD}^{cs}_{f^{-n}(x)})
$$

The trapping property $f(\overline{\cD}^{cs}_x)\en \cD^{cs}_{f(x)}$ implies directly that the sets $A_x$ verify:

\begin{itemize}
\item $f(A_x) = A_{f(x)}$ for every $x\in M$. \item The set $A_x$
is a decreasing intersection of topological balls (it is a cellular set).
In particular, $A_x$ is  compact and connected.
\end{itemize}

We would like to prove that the sets $A_x$ constitute a partition
of $M$, so that we can quotient the dynamics. For this, the
following lemma is of great use.

\begin{lema}\label{LemaRelacionEquivalencia}
For a given $x \in M$ and every $y\in \cW^{cs}(x)$, there exists $k_y$ such that
$f^{k}(\overline{\cD}^{cs}_y) \en \cD^{cs}_{f^{k}(x)}$ for every $k\geq k_y$. The
number $k_y$ can be chosen to vary semicontinuously on the point, that is, for every $y\in \cW^{cs}(x)$ there
exists $U$ a small neighborhood of $y$ (relative to $\cW^{cs}(x)$) such that for every $z\in
U$ we have that $k_z\leq k_y$.
\end{lema}

\dem{} Consider in $\cW^{cs}(x)$ the sets

$$ E_k = \{ y\in \cW^{cs}(x) \ : \ f^{k}(\overline{\cD}^{cs}_y) \en \cD^{cs}_{f^{k}(x)} \}$$

The sets $E_k$ are clearly open (by continuity of $f$ and
$\cD^{cs}$) and verify that $E_{k} \en E_{k+1}$ because of the
trapping property. Of course, $x \in E_k$ for every $k\geq 1$.

Using the trapping property and continuity of $\cD^{cs}$ again, it
follows that $\bigcup_{k\geq 0} E_k$ is closed. Indeed, let $y$
belong to the closure of $\bigcup_{k\geq 0} E_k$, so, for $z \in
\bigcup_{k\geq 0}E_k$ close enough to $y$ one has that
$f(\overline{\cD}^{cs}_y) \en \cD^{cs}_{f(z)}$. Since $z \in E_n$
for some $n$ one deduces that $y \in E_{n+1}$ showing that
$\bigcup_k E_k$ is closed. Being non-empty, one deduces that
$\bigcup_k E_k = \cW^{cs}(x)$ as desired.

The fact that the numbers $k_y$ can be chosen to vary semicontinuosly is a
consequence of the fact that $E_k$ is open ($k_y$ is the first
integer such that $y \in E_k$). \lqqd

In fact, a similar argument and continuity of the plaques gives a uniform estimate in the whole manifold.

\begin{lema}\label{l.uniformestimate} For every $n_0>0$ there exists $n_1>0$ such that for every $z\in
M$ and $w\in f^{n_0}\big(\cD^{cs}_{f^{-n_0}(z)}\big)$ we have that if $n\geq n_1$ then $$f^{n}(\cD^{cs}_{w}) \en \cD^{cs}_{f^{n}(z)}.$$
\end{lema}

\dem{} Notice that for each $z\in M$ there exists $n_z$ with
this property thanks to Lemma \ref{LemaRelacionEquivalencia} and the fact that
$D^{n_0}_z= f\big(\overline{\cD}^{cs}_{f^{-n_0}(x)}\big)$ is compact (it is enough to apply Lemma \ref{LemaRelacionEquivalencia} to each point of $D^{n_0}_z$ and use semicontinuity to get a uniform value of $n_z$ which works for every point in $D^{n_0}_z$). 

Now, we must show that the numbers $n_z$ vary semicontinuously with the point $z$. This follows from the fact that the disks $\cD^{cs}_w$ vary continuously: for $z'$ sufficiently close to $z$ and $w,w'$ in $f\big(\cD^{cs}_{f^{-n_0}(z)}\big)$, $f\big(\cD^{cs}_{f^{-n_0}(z')}\big)$ respectively, the disks $\cD^{cs}_{w}$ and $\cD^{cs}_{w}$ are very close. Therefore, their $n_z$'s iterates are close too and using continuity again one gets the desired property.  

Compactness of $M$ gives the existence of the desired $n_1$. 
\lqqd

One can now show that the sets $A_x$ constitute a partition of
$M$.

\begin{prop}\label{CorRelacionEquivalencia} For $x,y\in M$ we have that $A_x = A_y$ or $A_x \cap A_y =
\emptyset$.
\end{prop}

\dem{} Let $x \in M$ and consider $z \in A_x$.  

Since $z \in A_x \en \cW^{cs}(x)$, by Lemma \ref{LemaRelacionEquivalencia}, there exists $k$ such that $x \in f^{-k}\big(\cD^{cs}_{f^{k}(z)}\big)$. %Using the fact that $f(A_x)=A_{f(x)}$ we can assume that $x \in \cD^{cs}_z$.  

Lemma \ref{l.uniformestimate} gives $n_1>0$ such that given $n>0$ one has that 
$$f^{n_1}\big(\cD^{cs}_{f^{-n_1}(f^{-n+k}(x))}\big) \en
\cD^{cs}_{f^{-n+k}(z)}$$ which implies $A_{f^k(x)} \en
f^n(\cD^{cs}_{f^{-n+k}(z)})$ for every $n>0$. We have shown that if $z \in A_x$ then $A_{f^k(x)} \en A_{f^k(z)}$ for some $k>0$. Using the fact that $f(A_w)=A_{f(w)}$ for every $w$ we deduce that $A_x \en A_z$.    

This concludes because if $A_x \cap A_y \neq \emptyset$ this implies that there exists $z\in A_x \cap A_y$ and therefore $A_x \cup A_y \en A_z$. This gives that $x,y \in A_z$ and therefore, inverting roles one gets $A_z \en A_x \cap A_y$. We have shown $A_x \cup A_y \subset A_z \subset A_x \cap A_y$ which implies $A_x=A_y$ as desired.  
\lqqd

Consider two points $x,y$ such that $y \in \cW^{u}(x)$. We denote $\Pi^{uu}_{x,y}: \cD \en \cD^{cs}_x \to \cD^{cs}_y$ as the
unstable holonomy from a subset of $\cD^{cs}_x $ into a subset of $\cD^{cs}_y$. By the continuity of $\cD^{cs}$ and the bundle $E^u$ it holds that for $y$ close enough to $x$ the domain of $\Pi^{uu}_{x,y}$ is arbitrarily large. An important useful property is the following.

\begin{lema}\label{LemaHolonomia} The unstable holonomy preserves fibers, that is $\Pi^{uu}_{x,y}(A_x) = A_y$. In particular, $A_x$ is in the domain of $\Pi^{uu}_{x,y}$ for every $y$.
\end{lema}

\dem{}  It is enough to show (by the symmetry of the problem) that
$\Pi^{uu}(A_x)\en A_y$. For $n$ large enough we have that
$f^{-n}(\Pi^{uu}(A_x))$ is very close to a compact subset of
$\cD^{cs}_{f^{-n}(x)}$ and thus $f^{-n}(\Pi^{uu}(A_x)) \en
\cD^{cs}_{f^{-n}(y)}$ which concludes.

\lqqd

\begin{lema}\label{LemaSemicontinuidad}
The equivalence classes vary semicontinuously, i.e. if $x_n \to x$
then:
$$ \limsup A_{x_n} = \bigcap_{k>0} \overline{\bigcup_{n>k }
A_{x_n}} \en A_x $$
\end{lema}

\dem Using the invariance under unstable holonomy, it is enough to
show that the classes vary semicontinuously inside center-stable
manifolds. 

Consider $x\in M$ and $U$ a neighborhood of $A_x$ inside $\cW^{cs}(x)$. There exists $n>0$ such that $f^n\big(\overline{\cD}^{cs}_{f^{-n}(x)} \big) \en U$. Now, using Lemma \ref{l.uniformestimate} one deduces that if $w \in f^{n+1}\big(\cD^{cs}_{f^{-n-1}(x)} \big)$ then for some $k>0$ one has that $A_{f^k(w)} \en f^k(U)$ and therefore $A_w \en U$ as desired. 

\lqqd

We get a continuous projection by considering the relation
$x\sim y \Leftrightarrow y \in A_x$.

$$\pi : M \to M /_{\sim} $$

The space $M/_{\sim}$ with the quotient topology is Hausdorff thanks to the semicontinuity of the atoms of the partition (Lemma \ref{LemaSemicontinuidad}). In particular, since it is compact and second countable, it is metrizable.   

We denote as $g:M /_{\sim} \to M /_{\sim}$ the map given by
$g([x]) = [f(x)]$, that is 

$$g\circ \pi = \pi \circ f.$$ 

Since $\pi$ is continuous and surjective, it is a semiconjugacy.

Notice that a priori, the only knowledge one has on the topology of $M /_{\sim}$ is that it is the image by a cellular map of a manifold (some information on these maps can be found in the book \cite{Daverman} and references therein). For instance, we do not know a priori if the dimension of $M /_{\sim}$ is finite. This will follow from dynamical arguments after we prove Theorem \ref{TeoremaGesExpansivo} (combined with \cite{ManheExpansiveTopological}).

Given a homeomorphism $h: X \to X$ of a compact metric space $X$, we denote the $\eps$-stable ($\eps$-unstable) set as

$$ S_\eps  (x)= \{ y \in X \ : \  d(h^n(x), h^n(y)) \leq \eps \ for \  n \geq 0 \}$$

$$ U_\eps (x)=  \{ y \in X \ : \  d(h^n(x), h^n(y)) \leq \eps \ for \  n \leq 0 \}$$

We say that a homeomorphism has \emph{local product structure} if there exists $\delta>0$ such that $d(x,y)<\delta$ implies that $S_\eps (x)\cap U_\eps (y) \neq \emptyset$.

A homeomorphism $h$ is \emph{expansive} (with expansivity constant $\alpha>0$) if for every $x\in X$ we have that $S_\alpha (x) \cap U_\alpha(x) = \{x\}$.

Expansive homeomorphisms verify that $\diametro(h^n(S_\eps(x))) \to 0$ uniformly on $x$ for $\eps<\alpha$.

We shall also denote $W^{s(u)}(x) = \{ y\in X \ : \ d(h^n(x),h^n(y)) \to 0 \ as \ n\to +\infty (n\to -\infty) \}$, the above remark implies that $S_\eps(x) \en W^s(x)$ for an expansive homeomorphism ($\eps < \alpha$).

It is important to remark that expansivity is a purely topological notion and is independent of the chosen metric on $X$. The following is classical.

\begin{lema}\label{l.expansiveequivalence} A homeomorphism $h$ of a compact metric space $X$ is expansive with local product structure if and only if 
there exists two neighborhoods $V_1\en V_2$ of the diagonal $\Delta \en X \times X$ with the following properties: 
\begin{itemize}
\item the maximal invariant set of $h\times h$ in $V_2$ is $X\times X$, i.e. if $(x,y) \in V_2$ is such that $(h^n(x),h^n(y)) \in V_2$ for every $n\in \ZZ$ then $x=y$. 
\item if $(x,y) \in V_1$ then, if $S_{V_2}(x)= \{ z \in X \ : \ (h^n(x),h^n(z)) \in V_2 \ , \ \forall n \geq 0 \}$ and $U_{V_2}(x)=  \{ z \in X \ : \ (h^n(z), h^n(y)) \in V_2 \ , \ \forall n \leq 0 \}$ one has: 
$$   S_{V_2}(x) \cap U_{V_2}(y) \neq \emptyset $$  
\end{itemize}
\end{lema}

Notice that the first condition is equivalent to $S_{V_2}(x) \cap U_{V_2}(x)= \{x\}$. 

The proof is left to the reader, it is a direct consequence of compactness of $X$. See also \cite{Ombach} and references therein for properties of such homeomorphisms which are sometimes called \emph{Smale homeomorphisms} (they behave topologically exactly as hyperbolic diffeomorphisms). 

\begin{teo}\label{TeoremaGesExpansivo}
The homeomorphism $g$ is expansive with local product structure.
Moreover, $\pi(\cW^{cs}(x))= W^s(\pi(x))$ and $\pi$ is injective
when restricted to the unstable manifold of any point.
\end{teo}

\dem The last two claims are direct from Lemma
\ref{LemaRelacionEquivalencia}, Lemma \ref{LemaHolonomia} and the definition of the
equivalence classes.

We must show the existence of a local product structure and that will establish expansivity also. First choose $\eps>0$ such that an unstable manifold of size $2\eps$ cannot intersect the same center stable disk in more than one point. This is given by the continuity of the bundles $E^{cs}$ and $E^u$.

Consider $x \in M$ and a neighborhood $U$ of $A_x$. Using Lemma \ref{LemaSemicontinuidad} one knows that there is a neighborhood $V$ of $A_x$ such that for every $y \in V$ one has that $A_y$ is contained in $U$. One can choose $U$ small enough so that for every $y \in V$ one has that $\cW^{u}_{\eps}(y) \cap \cD^{cs}_x$ is exactly one point.

Moreover, by the continuous variation of the $\cD^{cs}$-disks, one has that, maybe by choosing $V$ smaller, it holds that for every $y,z\in V$ one has that $\cW^{u}_{2\eps}(y) \cap \cD^{cs}_z$ is exactly one point.

Since the image of $V$ by $\pi$ is open, one gets a covering of $M/_\sim$ by open sets where there is product structure in the sense of Lemma  \ref{l.expansiveequivalence}. By compactness one deduces that there exists a local product structure for $g$ and since the intersection point is unique one
also obtains expansivity of $g$. 

\lqqd

\subsection{Some remarks on the topology of the quotient} This
section shall not be used in the remaining of the paper.

We shall cite some results from \cite{Daverman} which help
understand the topology of $M /_{\sim}$. Before, we remark that
Ma\~n\'e proved that a compact metric space admitting an expansive
homeomorphism must have finite topological dimension
(\cite{ManheExpansiveTopological}).

Corollary IV.20.3A of \cite{Daverman} implies that, since $M
/_{\sim}$ is finite dimensional, we have that it is a locally
compact ANR (i.e. absolute neighborhood retract). In particular,
we get that $\dim (M/_{\sim}) \leq \dim M$ (see Theorem III.17.7).
Also, using Proposition VI.26.1 (or Corollary VI.26.1A) we get
that $M/_{\sim}$ is a $d-$dimensional homology manifold (since it
is an ANR, it is a \emph{generalized manifold}). More properties
of these spaces can be found in section VI.26 of \cite{Daverman}.

Also, in the cited book, one can find a statement of Moore's
theorem (see section IV.25 of \cite{Daverman}) which states that a
cellular decomposition of a surface is approximated by
homeomorphisms. In particular, in our case, if $\dim E^{cs} = 2$,
we get that $M /_{\sim}$ is a manifold (see also Theorem VI.31.5
and its Corolaries).

Later, we shall see that if $M=\TT^d$ then the quotient $M/_\sim$
is also a manifold (indeed $M/_\sim$ is homeomorphic to $\TT^d$).
The same should hold for infranilmanifolds but we have not checked
this.
%
%Some other results are available, in particular, we notice
%Edward's cell-like decomposition theorem which states that if
%$\sim$ is a cellular decomposition of a $d$ dimensional manifold
%($d\geq 5$) such that $M /_\sim$ has finite topological dimension
%and such that it has the \emph{disjoint disk property} (see
%chapter IV.24 of \cite{Daverman}) then the quotient map is
%aproximated by homeomorphisms. A similar result exists for
%dimension $3$ which is even more technical. Notice than in our
%case, since we have the decomposition of the center-stable
%manifold, we can play with the dimensions in order not to be never
%in dimension $4$ by choosing to work with the decomposition on the
%center stables or the whole manifold.
%
%Also, we remark that it is known that when multiplying a
%decomposition by $\RR^2$ we always get the \emph{disjoint disc
%property} and in all known decompositions, after multiplying by
%$\RR$ we get a decomposition approximated by homeomorphisms (see
%section V.26 of \cite{Daverman}) so in theory, it should be true
%that always our space $M  /_\sim$ is a manifold homeomorphic to
%$M$. 

\section{Transitivity of the expansive
homeomorphism}\label{Section-TransitiveExpansive}

In this section $g: M /_\sim \to M/_\sim$ will denote the expansive quotient map we have constructed in the previous section. The quotient map will be as before denoted by $\pi: M \to M/_\sim$.

It is not yet known if an Anosov diffeomorphism must be transitive. Since Anosov diffeomorphisms enter in our hypothesis, there is no hope of knowing if $f$ or $g$ will be transitive without solving this long-standing conjecture. We shall then work with similar hypothesis to the well known facts for Anosov diffeomorphisms, showing that those hypothesis that we know guarantee that Anosov diffeomorphisms are transitive imply transitivity of $g$ as defined above.

\begin{obs} It is well known that transitivity of $g$ amounts to showing some form of uniqueness of basic pieces. This is quite direct if one assumes some knowledge on the structure of the foliations of $f$, for example, if for every $x,y \in M$ one has that $\cD^{cs}_x \cap \cW^u(y) \neq \emptyset$ then it follows that $g$ is transitive. In this paper we rather concentrate on information which does not rely \emph{a priori} on knowledge of the structure of the foliations.
\end{obs}

In particular, we shall prove in this section the following two results.

\begin{teo}\label{TeoremaNewhouse} Assume $f: M \to M$ is a dynamically coherent partially hyperbolic diffeomorphism with the trapping property and $\dim E^u=1$. Then $M$ is covered by $\RR^d$ and homotopically equivalent to $\TT^d$.
\end{teo}

For $d=3$ or $d\geq 5$ a manifold homotopically equivalent to $\TT^d$ is indeed homeomorphic to $\TT^d$ (by geometrization for $d=3$ and see \cite{HsiangWall} for $d\geq 5$). Even if there are counterexamples for this purely topological result for $d=4$ this does not matter since dynamical arguments will give us that the manifold must be topologically $\TT^d$ anyway. 

\begin{teo}\label{TeoremaGtransitivo} Assume $f: M \to M$ is a dynamically coherent partially hyperbolic diffeomorphism with the trapping property and $M$ is covered by $\RR^d$ and homotopically equivalent to $\TT^d$. Then $f$ is isotopic to a linear Anosov automorphism $L$, the manifold $M$ and the quotient $M/_\sim$ are homeomorphic to $\TT^d$ and $g$ is topologically conjugate to $L$.
\end{teo}

Put together, Theorems \ref{TeoremaNewhouse} and \ref{TeoremaGtransitivo} can be compared to Franks-Newhouse theory (\cite{Franks,Newhouse}) on codimension one Anosov diffeomorphisms. It is possible to prove directly (with an argument similar to the one of Newhouse but taking care on the quotients) that $g$ is transitive when $\dim E^u=1$ without showing that $M = \TT^d$ (see \cite[Section 5.4]{Pot-Tesis} for this approach).

Theorem \ref{TeoremaGtransitivo} is reminiscent of Franks-Manning theory (\cite{Franks,Manning} see also \cite[Chapter 18.6]{KH}). It is natural to expect that property this result should hold if we consider $M$ an infranilmanifold, but we have not checked this in detail.

It is reasonable to extend the conjecture about transitivity of Anosov diffeomorphisms to expansive homeomorphisms in manifolds with local product structure. See the results in \cite{Vieitez,ABP}.

\subsection{Proof of Theorem \ref{TeoremaNewhouse}}

This proof is an adaptation of quite classical ideas (see for example the Appendix in \cite{ABP}) with some arguments of \cite{Newhouse}.

One key point is that $\cW^{cs}$ is a foliation by leafs homeomorphic to $\RR^{d-1}$ which follows directly from the trapping property (and Lemma \ref{LemaRelacionEquivalencia}) giving that the leafs of $\cW^{cs}$ are increasing union of disks. 

Having this, one can lift the foliations $\cW^{cs}$ and $\cW^{u}$ to the universal cover and show that a leaf of $\tilde \cW^{u}$ cannot
intersect a leaf of $\tilde \cW^{cs}$ more than once using Haefliger's argument (\cite[Proposition 7.3.2]{CandelConlon}) and the fact that all
leafs are simply connected.

To prove that the universal cover of $M$ is $\RR^d$ one must show that given a leaf of $\tilde \cW^{u}$ it intersects every leaf of $\tilde \cW^{cs}$. This follows with exactly the same proof of Lemma (5.2) of \cite{Franks} once one knows that every leaf of $\cW^{cs}$ is dense.

\begin{lema}\label{Lemma-LeafsAreDense} Every leaf of $\cW^{cs}$ is dense in $M$.
\end{lema}

\dem We use here some of the ideas of \cite{Newhouse}.

Consider a spectral decomposition for $g$ (the proof for hyperbolic diffeomorphisms can be seen in \cite{NewhouseHyperbolic}, the same proof applies for expansive homeomorphisms with local product structure, see again \cite{Ombach} and references therein). In this way, one sees that the foliation $\cW^{cs}$ contains finitely many minimal sets (associated to the basic pieces of $g$ which are repellers).

Consider a set $\Lambda$ which is the preimage by $\pi$ of a repeller of $g$. The set $\Lambda$ is $f$-invariant and consists of finitely many minimal sets of the foliation $\cW^{cs}$. It is enough to show that $\Lambda=M$, for this, it is enough to show that every point $x\in \Lambda$ verifies that $\cW^{u}(x) \cap \Lambda \neq \emptyset$ in both connected components of $\cW^{u}(x)\setminus \{x\}$ (recall that since $\cW^{u}(x)$ is one-dimensional $\cW^{u}(x) \setminus \{x\}$ has two different connected components). Under this assumption, and using the fact that $f^{-1}$ contracts uniformly leafs of $\cW^u$ one deduces that $\Lambda$ is also saturated by $\cW^u$ showing that $\Lambda$ is open and closed. Being non-empty, one deduces that $\Lambda=M$ as desired.

So, to prove the lemma it is enough to show the following.

\begin{af} Every point $x \in  \Lambda$ verifies that both connected components of $\cW^u(x)\setminus \{x\}$ intersect $\Lambda$.
\end{af}

\dem We claim that if $x \in \Lambda$ is a point such that one connected component of $\cW^{u}(x) \setminus \{x\}$ does not intersect $\Lambda$ then $A_x$ is periodic (or equivalently, $\pi(x)$ is a periodic point for $g$). Moreover, there are finitely many such periodic points. To see this, notice that there
exists $\eps>0$ such that if three points of the past orbit of $\pi(x)$ by $g$ are at distance smaller than $\eps$ one has that the unstable manifold of one of the points in the orbit intersects $\Lambda$ in both connected components. Since such points are invariant one deduces that $\pi(x)$ must be periodic for $g$.

Now assume there is a point $x$ such that its unstable manifold does not intersect $\Lambda$ on one side. Let $\Sigma$ be the boundary of $\cD^{cs}_x$ which is a topological sphere (of dimension $\geq 1$ since $d\geq 3$). As we mentioned, every point in $\Sigma$ verifies that the unstable manifold in both sides
intersect $\Lambda$, and by continuity and the intersection point, one obtains a continuous map from $\varphi: \Sigma \times [0,1] \to M$ which verifies that $\varphi(z,0) = z$, $\varphi(z, 1)$ is in $\Lambda$ and maps $\{z\} \times [0,1]$ to a compact part of $\cW^u(z)$. We can moreover assume that $\varphi(\Sigma \times \{t_0\}$ is contained in a leaf of $\cW^{cs}$ for every $t_0$ using continuity (see \cite{Newhouse}).

Since $\varphi(\Sigma \times \{t_0\})$ separates $\cW^{cs}(\varphi(z,t_0))$ giving a compact region one can prove that the unstable manifold of $x$ intersects
$\cW^{cs}(\varphi(z,t_0))$ for every $t_0 \in [0,1]$ and therefore that it intersects $\Lambda$ giving a contradiction. See \cite{Newhouse} for more details.
\lqqd

Now, we have a global product structure in the universal cover which implies that $\tilde M = \RR^d$ and moreover, we get that the space of leafs of the foliation $\tilde \cW^{cs}$ is homeomorphic to the real line $\RR$ (and can be identified with a single leaf of $\tilde \cW^u$). The action by deck transformations induces an action on the space of leafs of $\tilde \cW^{cs}$ which does not have fixed points since all leafs of $\cW^{cs}$ are
simply connected. By H\"{o}lder's theorem this implies that $\pi_1(M)$ is free abelian\footnote{H\"{o}lder's theorem implies that the group is semiconjugate to a subgroup of translations of the real line. In particular, it has no torsion. Since $\pi_1(M)$ is finitely presented because $M$ is compact, one deduces that it is $\ZZ^k$ for some $k\geq 0$.} and thus isomorphic to $\ZZ^k$. Since the universal cover of $M$ is contractible, it is a $K(\ZZ^k,1)$ and therefore, it is homotopy equivalent to $\TT^k$ and since $M$ is a compact manifold, $k=d$ and $M$ is homotopy equivalent to $\TT^d$ as desired.

\lqqd

\subsection{Proof of Theorem \ref{TeoremaGtransitivo}} We shall follow the proof given in
\cite{KH} chapter 18.6.

Before we start with the proof, we shall recall Theorem 18.5.5 of \cite{KH} (the statement is modified in order to fit our needs,
notice that for an expansive homeomorphism with local product structure the specification property is verified in each basic
piece, see \cite{Ombach} and references therein):

\begin{prop}[Theorem 18.5.5 of \cite{KH}]\label{PropConteoDePeriodicos}
Let $X$ a compact metric space and $g:X \to X$ an expansive homeomorphism with local product structure. Then, there exists
$h,c_1,c_2>0$ such that for $n\in \NN$ we have:
$$ c_1 e^{n h} \leq P_n(g) \leq c_2 e^{n h} $$
\noindent where $P_n(g)$ is the number of fixed points of $g^n$.
\end{prop}

We shall use several time the very well know Lefschetz formula which relates the homotopy type of a continuous function, with the
index of its fixed points (see \cite{FranksLibro} Chapter 5).

\begin{defi}
Let $V\en \RR^k$ be an open set, and $F: V\en \RR^k \to \RR^k$ a continuous map such that $\Gamma\en V$ the set of fixed points of
$F$ is a compact set, then, $I_\Gamma(F) \in \ZZ$ (the index of $F$) is defined to be the image by $(id-F)_\ast: H_k(V, V-\Gamma)
\to H_k(\RR^k, \RR^k-\{0\})$ of $u_\Gamma$ where $u_\Gamma$ is the image of $1$ under the composite $H_k(\RR^k , \RR^k - D) \to
H_k(\RR^k, \RR^k- \Gamma) \cong H_k(V, V-\Gamma)$ where $D$ is a disk containing $\Gamma$.
\end{defi}

\begin{obs}
In general, if we have a map $h$ from a $d$-dimensional manifold $M$ to itself, we can embed the manifold in $\RR^k$ by a map $\imath : M \hookrightarrow \RR^k$ for some big $k>d$ and one has a retraction $r: V \to M$ of a neighborhood $V$ of $\imath(M)$. The value of $I_\Gamma(f) = I_\Gamma(\imath \circ f \circ r)$ does not depend on the embedding nor the retraction. This is also equivalent to consider the fixed point set in a chart of $M$ and computing the index in the chart. 
\end{obs}

One can also see that if $\Gamma= \mathrm{Fix}(h) = \Gamma_1 \cup \ldots \cup \Gamma_k$ where $\Gamma_i$ are compact and disjoint, then $I_\Gamma(h)= \sum_{i=1}^k I_{\Gamma_i}(h)$. Here we consider $I_{\Gamma_i}(h)$ as the index restricted to an open set $V_i$ of $\Gamma_i$ which does not intersect the rest of the $\Gamma_j$. See \cite[Theorem 5.8 (b)]{FranksLibro}.

For a single hyperbolic fixed point, it is very easy to compute the index, it is exactly $sgn(det(Id - D_pf))$ (\cite[Proposition 5.7]{FranksLibro}). Since the definition is topological, any time we have a set which behaves locally as a hyperbolic fixed point, it is not hard to see that the index is
the same (\cite[Theorem 5.8 (c)]{FranksLibro}).

Lefschetz fixed point formula (\cite[pages 34-38]{FranksLibro}) for the torus can be stated as follows:

\begin{teo}[Lefschetz fixed point formula]\label{teoLefshetz}
Let $h:\TT^d \to \TT^d$ be a homeomorphism with fixed point set $\Gamma=\mathrm{Fix}(h)$. Then, the index $I_\Gamma(h)= det(Id- h_\ast)$ where $h_\ast: H_1(\TT^d,\ZZ) \to H_1(\TT^d, \ZZ)$ is the action of $h$ in homology.
\end{teo}

Now we come back to the proof of Theorem \ref{TeoremaGtransitivo}. The first thing we must show, is that the \emph{linear} part of
$f$, that is, the action $L= f_\ast: H_1(\TT^d, \ZZ) \to H_1(\TT^d, \ZZ) \in SL(d,\ZZ)$ is a hyperbolic matrix.

\begin{lema}\label{LemaParteLinealHiperbolica} The matrix $L$ is hyperbolic.
\end{lema}

\dem We can assume (maybe after considering a double covering and $f^2$) that $E^{cs}$ and $E^u$ are orientable and its orientations preserved by $Df$. So, it is not hard to show that for every fixed point $p$ of $g^n$, the index of $\pi^{-1}(p)$ for $f$ is of modulus one and always of the same sign.

So, we know from the Lefshetz formula that

$$|det(Id - L^n)| = \sum_{g^n(p)=p} |I_{\pi^{-1}(p)}(f)| = \# \mathrm{Fix}(g^n).$$

This implies that $L^n$ is hyperbolic using Proposition
\ref{PropConteoDePeriodicos} since the only way to have that
estimate on the periodic orbits is that $L$ is hyperbolic (see the
argument in Lemma 18.6.2 of \cite{KH}). \lqqd

It is standard to show the existence of a semiconjugacy $h:\TT^d
\to \TT^d$ isotopic to the identity such that $h \circ f = L \circ
h$. Its lift $\overline h: \RR^d \to \RR^d$ is given by shadowing,
in particular, the iterations of the set $(\overline h)^{-1}(x)$
remain of bounded diameter.

\begin{lema}\label{LemaSemiconjIntermedia}
The semiconjugacy $h$ factors through the quotient map $\pi$.
More precisely, there exists $\tilde h : \TT^d /_{ \sim} \to
\TT^d$ continuous such that $\tilde h \circ \pi = h$.
\end{lema}

\dem It is enough to show that for every $x\in \TT^d/_\sim$ there
exists $y\in \TT^d$ such that $\pi^{-1}(x) \subset h^{-1}(y)$.

For this, notice that any lifting  of $\pi^{-1}(x)$ (that is, a
connected component of the preimage under the covering map) to the
universal covering $\RR^d$ verifies that it's iterates remain of
bounded size. This concludes by the remark above on $\overline h$.
\lqqd

Now, we shall prove that if $\overline f : \RR^d \to \RR^d$ is any
lift of $f$, then there is exactly one fixed fiber of $\pi$ for
$\overline f$.

\begin{lema}\label{LemaUnicaFibraArriba}
Let $\overline f^n$ be any lift of $f^n$ to $\RR^d$. So, there is
exactly one fixed fiber of $\pi$.
\end{lema}

\dem Since $\overline f^n$ is homotopic to $L^n$ which has exactly
one fixed point and each fixed fiber of $\pi$ contributes the same
amount to the index of $\overline f^n$ it must have exactly one
fixed fiber. \lqqd

This allows us to show that $g$ is transitive:

\begin{prop}\label{PropGtransitivo} The homeomorphism $g$ is transitive.
\end{prop}

\dem First, we show that there exists a basic piece of $g$ which
projects by $\tilde h$ to the whole $\TT^d$. This is easy since
otherwise, there would be a periodic point $q$ in $\TT^d
\backslash \tilde h (\Omega(g))$ but clearly, the $g-$orbit of
$\tilde h^{-1}(q)$ must contain non-wandering points (it is
compact and invariant). This concludes, since considering a
transitive point $y$ of $L$ and a point in $\Omega(g) \cap \tilde
h^{-1}(y)$ we get the desired basic piece.

Now, let $\Lambda$ be the basic piece of $g$ such that $\tilde h
(\Lambda) = \TT^d$. Assume that there exists $\tilde \Lambda \neq
\Lambda$ a different basic piece and $z$ a periodic point of
$\tilde \Lambda$, naturally, we get that $\tilde h^{-1}(\tilde h
(z))$ contains also a periodic point $z'$ in $\Lambda$. By
considering an iterate, we can assume that $z$ and $z'$ are fixed
by $g$.

So, we get that lifting $h^{-1}(\tilde h(z))$ to a lift which
fixes $\pi^{-1}(z)$ and $\pi^{-1}(z')$ which contradicts the
previous lemma.

\lqqd

With this in hand, we will continue to prove that the fibers of
$h$ coincide with those of $\pi$ proving that $g$ is conjugate to
$L$ (in particular, $\TT^d /_\sim \cong \TT^d$).

First, we show a global product structure for the lift of $f$.
Notice that when we lift $f$ to $\RR^d$, we can also lift its
center-stable and unstable foliation. It is clear that both
foliations in $\RR^d$ are composed by leaves homeomorphic to
$\RR^{cs}$ and $\RR^u$ respectively (the unstable one is direct,
the other is an increasing union of balls, so the same holds).

\begin{lema}\label{LemaGlobalProductStructure}
Given $x,y\in \RR^d$, the center stable leaf of $x$ intersects the
unstable leaf of $y$ in exactly one point.
\end{lema}

\dem The fact that they intersect in at most one point is given by
the fact that otherwise, we could find a horseshoe for the lift,
and thus many periodic points contradicting  Lemma
\ref{LemaUnicaFibraArriba} (for more details, see Lemma 18.6.7 in
\cite{KH}).

The proof that any two points have intersecting manifolds, is
quite classical, and almost topological once we know that both
foliations project into minimal foliations (see also Lemma 18.6.7
of \cite{KH}). \lqqd

Now, we can conclude with the proof of Theorem
\ref{TeoremaGtransitivo}.

To do this, notice that the map $\overline h$ conjugating
$\overline f$ with $L$ is proper, so the preimage of compact sets
is compact. Now, assume that $A_x, A_y$ are lifts of fibers of
$\pi$ such that $\overline h (A_x) = \overline h (A_y)$ we shall
show they coincide.

Consider $K$ such that if two points have an iterate at distance
bigger than $K$ then their image by $\overline h$ is distinct.

We fix $x_0 \in A_x$ and consider a box $D_K^n$ of $\overline f^n
(x_0)$ consisting of the points $z$ of $\RR^d$ such that $\cW^u(z)
\cap \cW^{cs}_K(x_0) \neq \emptyset$ and $\cW^{cs}(z) \cap
\cW^u_K(x_0) \neq\emptyset$.

It is not hard to show using Lemma
\ref{LemaGlobalProductStructure} that there exists $\tilde K$
independent of $n$ such that every pair of points in $D_K^n$ in
the same unstable leaf of $\cW^u$ have distance along $\cW^u$
smaller than $\tilde K$ (this is a compactness argument). An
analogous property holds for $\cW^{cs}$.

This implies that if $\overline f^n (A_y) \en D_K^n$ for every
$n\in \ZZ$ then $A_y$ and $A_x$ must be contained in the same leaf
of $\cW^{cs}$. In fact we get that $(\overline{f})^{-n}(A_y) \en
\cW^{cs}_K((\overline{f})^{-n}(x_0))$  for every $n\geq 0$ and so we
conclude that $A_x=A_y$ using Lemma
\ref{LemaRelacionEquivalencia}.

\lqqd

\subsection{Some manifolds which do not admit this kind of diffeomorphisms}\label{SubSubsectionVARIEDADESQADMITEN}

The arguments used in the previous section also allow to show that
certain manifolds (and even some isotopy classes in some
manifolds) do not admit dynamically coherent partially hyperbolic
diffeomorphisms satisfying the trapping property.
%
%To do this, we recall that for general manifolds $M^d$, and a
%homeomorphism $h: M\to M$, the Lefschetz number of $f$, which we
%denote as $L(h)$ is calculated as $\sum_{i=0}^d trace(f_{\ast,i})$
%where $h_{\ast,i}: H_i(M,\QQ) \to H_i(M,\QQ)$ is the induced map
%on (rational) homology. We also have that the sum of index the sum
%of the Lefshetz index along a covering of $Fix(h)$ by sets
%homeomorphic to balls equals $L(h)$.

A similar argument to the one used in the previous section yields
the following result (see \cite{GogolevHertz} for sharper results
in the same lines).

\begin{teo}\label{TeoremaAlgunasVariedadesNoAdmiten}
Let $f$ be a partially hyperbolic diffeomorphism of $M$ with the
coherent trapping property, then, the action $f_\ast :
H_{\ast}(M,\QQ) \to H_{\ast}(M,\QQ)$ is strongly partially
hyperbolic (it has both eigenvalues of modulus $>1$ and $<1$).
\end{teo}

This leads to a natural question: is every dynamically coherent
partially hyperbolic diffeomorphism with the trapping property
homotopic to an Anosov diffeomorphism?. One should notice that
expansive homeomorphisms admitting transverse stable and unstable
foliations share many properties with Anosov diffeomorphisms (see
for example \cite{Vieitez,ABP}) but it is not known if every such
homeomorphism is topologically conjugate to an Anosov
diffeomorphism.

Also, let us remark that there exist examples of dynamically coherent
partially hyperbolic diffeomorphisms which are isotopic to Anosov
and robustly transitve while not satisfying the trapping property.
See \cite[Section 3.3.4]{Pot-Tesis}.

%%%%%%%%%%%%%%%%%%%%%%%%%%%%%%%%%%%%%%%%%%%%%%%%%%%%%%%%%%%%%%%%%%%%%%%%%%%%%%%%%%%%%%%%%%%%%%%%%%%%%%%%%%%%

\section{Some dynamical consequences}\label{Section-Dynamics}

In this section we shall look at what type of dynamical properties
can be recovered in the spirit of \cite{Pot-WildMilnor} (see also \cite{Carvalho,BV}).

We recall that a \emph{quasi-attractor} $\Lambda$ is a chain-recurrence class
satisfying that it admits a decreasing basis of neighborhoods $U_n$ satisfying
that $f(\overline{U_n})\en U_n$ (see \cite[Chapter 1]{Pot-Tesis} and references therein).

Since a quasi-attractor is saturated by unstable manifolds and the quotient we have defined which
conjugates $f$ to an expansive homeomorphism is injective on unstable manifolds, one expects that
whenever the quotient map $g$ is transitive (as it is ensured in some cases by Theorem \ref{teoPPal})
there is a unique quasi-attractor. Unfortunately, showing this would involve showing that
there are fibers of the semiconjugacy which are trivial and this is a subtle issue as the example presented
in Appendix \ref{Appendix-Decomposition} shows.

We are however able to show uniqueness of the quasi-attractors under a mild assumption
resembling \emph{chain-hyperbolicity} as defined in \cite{CP}.

\begin{prop}\label{Prop-UniqueQuasiAttractor}
Let $f: M \to M$ be a dynamically coherent partially hyperbolic diffeomorphism satisfying the trapping property.
Assume moreover that the quotient map $g$ defined above is transitive and that there exists a point $x\in M$ such
that $A_x=\{x\}$. Then, $f$ has a unique quasi-attractor.
\end{prop}

\dem Consider $\Lambda$ a quasi-attractor for $f$ and let $\pi: M \to M/_\sim$ be the semiconjugacy to $g: M/_\sim \to M/_\sim$ constructed
in section \ref{Section-ExpansiveQuotient}.

Since $\pi$ is injective along unstable manifolds, one obtains that $\pi(\Lambda)$ contains the unstable set of any point
$z \in M/_\sim$ such that $z = \pi(y)$ with $y\in \Lambda$. Since $g$ is expansive with local product structure and transitive, one know that the
orbit of $W^u(z)$ is dense in $M/_\sim$ (see for example \cite[Chapter 18]{KH}).

This implies that $\pi(\Lambda)$ is dense in $M/_\sim$. Since $\pi$ is continuous and $\Lambda$ compact one deduces that $\pi(\Lambda)= M/_\sim$. In particular, $\pi(\Lambda)$
contains $\pi(x)$.

As a consequence, we get that every quasi-attractor must intersect $A_x=\{x\}$. Since different quasi-attractors must be disjoint, this implies uniqueness of the quasi
attractor under the assumptions of the proposition.

\lqqd

\begin{obs}
 In the case where $E^{cs}= E^s\oplus E^c$ with $\dim E^{c}=1$ and a trapping property is verified by leaves tangent to $E^c$ one can show that $f$ satisfies a trapping property.
 From the construction of $\pi$ one sees that $A_x$ is either a point or a closed interval. Using \cite{Jones} one sees immediately that the conditions of the
 previous proposition are satisfied.
\end{obs}

Another property, related to \cite{Pot-WildMilnor} is the following.

\begin{prop}\label{Prop-ClasesEnSuperficies}
  Let $f: M \to M$ be a dynamically coherent partially hyperbolic diffeomorphism satisfying the trapping property.
  Assume moreover that the quotient map $g$ defined above is transitive and the image of every open set of a center-stable leaf by $\pi$
  is either a point or has non-empty interior in the stable manifold of $g$. Then $f$ has a unique quasi-attractor and every other
  chain-recurrence class of $f$ is contained in a periodic disk of $\cW^{cs}$.
\end{prop}

\dem Consider a quasi-attractor $\Lambda$. As in the previous proposition, one has that $\pi(\Lambda)= M/_\sim$ and it is saturated by $\cW^u$.

One can easily show that for every $x \in M$ the boundary of $A_x$ is contained in $\Lambda$: indeed, consider any $y \in \partial A_x$ and
a neighborhood $U$ of $y$ in $\cW^{cs}(y)$. From our hypothesis one has that $\pi(U)$ has non-empty interior in the stable manifold of $\pi(y)$.

Iterating backwards and using the semiconjugacy and using the density of unstable sets for $g$ one obtains that $f^{-n}(U)$ intersects $\Lambda$ for some
$n$. Invariance of $\Lambda$ and the fact that the choice of $U$ was arbitrary gives that $y \in \overline{\Lambda}=\Lambda$.

The rest of the proposition follows by applying Proposition 2.1 of \cite{Pot-WildMilnor}.

\lqqd

As we have explained, the hypothesis we demand in this section might follow directly from the fact that $f$ has the trapping property but
the example presented in Appendix \ref{Appendix-Decomposition} strongly suggests that counterexamples might exist.

\begin{quest}
 Does there exists a dynamically coherent partially hyperbolic diffeomorphism of $\TT^3$ with splitting $T\TT^3 = E^{cs} \oplus E^u$ and with the trapping property such that it admits more than one quasi-attractor?
 Such that it has chain-recurrence classes (different than the quasi-attractor) which are not contained in periodic center-stable discs?
\end{quest}

See \cite{Pot-JDDE} for related discussions.

%%%%%%%%%%%%%%%%%%%%%%%%%%%%%%%%%%%%%%%%%%%%%%%%%%%%%%%%%%%%%%%%%%%%%%%%%%%%%%%%%%%%%%%

\section{A weaker trapping property and
coherence}\label{Section-Coherence}

In this section we shall present a weaker trapping property
without requiring dynamical coherence \emph{a priori} and show
that it is enough to recover the initial proposition. One would
hope that this property is shared by certain partially hyperbolic
diffeomorphisms isotopic to Anosov though it is not so clear that
it holds (see \cite{Pot-PHFoliations} for results in this
direction). The proof is completely analogous to the one presented
in section 3 of \cite{BF} but in a slightly different context. One
important point is the fact that we do not assume that the
trapping property occurs in a small region and hope this might
find applications.

Let $f: M \to M$ be a partially hyperbolic diffeomorphism with
splitting $TM = E^{cs} \oplus E^u$.

As before, we denote as $cs=\dim E^{cs}$ and $u= \dim E^u$ and
$\DD^{\sigma}$ is the $\sigma$-dimensional open disk and
$\overline{\DD}^{\sigma}$ its closure.

We will assume that $f$ verifies the following property:

\begin{itemize}
\item[($\ast$)] there exists a continuous map $B: M \to
\Emb^1(\overline{\DD}^{cs} \times \overline{\DD}^u)$ such that
$B(x)(0,0)=x$ for every $x \in M$, one has that $B(x)( \{a\}
\times \overline{\DD}^u) \en \cW^u(B(x)(a,0))$ and the following
trapping property is verified:

$$ \forall a \in \overline{\DD}^{cs} \ \ \exists b \in \DD^{cs} \
\ ; \ \ B(f(x))(\{b\} \times \overline{\DD}^u) \en f(B(x)(\{a\}
\times \DD^u)) .$$
\end{itemize}

The main result of this section is the following.

\begin{teo}\label{Teo-TrappingImplicaCoherencia} Let $f: M \to M$
be a partially hyperbolic diffeomorphism verifying property
$(\ast)$, then, $f$ is dynamically coherent with a trapping
property.
\end{teo}

\dem First, we will denote as $\cD^{cs}_x$ to the set of points $y \in B(x)(\DD^{cs} \times \DD^u)$ such that

$$ f^n(y) \in B(f^n(x))(\DD^{cs} \times \DD^u) $$

We claim that $\cD^{cs}_x$ is a manifold everywhere tangent to
$E^{cs}$ and moreover one has that a trapping property
$f(\overline{\cD}^{cs}_x) \en \cD^{cs}_{f(x)}$ is verified. Also,
$\cD^{cs}_x$ intersects every local unstable manifold of $\Image
B_x= B(x)(\DD^{cs} \times \DD^u)$ in a unique point.

To show this, notice first that expansivity of the unstable
manifolds implies that $\cD^{cs}_x$ cannot intersect an unstable
manifold more than once. Also, the trapping property verified by
the maps $B$ gives that every point $a \in \DD^{cs}$ verifies that
$B(x)(\{a\}\times \DD^u)$ intersects $\cD^{cs}_x$. The fact that
it is a $C^{1}$ manifold everywhere tangent to $E^{cs}$ follows by
classical graph transform arguments (see \cite{HPS} or
\cite[Chapter 6]{KH}).

An important fact of the above is that one can view $\cD^{cs}_x$
as a limit of disks $D^n_x$ where $D^n_x$ is any disk inside
$\Image B_x$ with the following property.

\begin{itemize}
\item There exists a disk $D$ in $\Image B_{f^n(x)}$ which is the
image by $B(f^n(x))$ of a $1$-Lipchitz graph over
$\DD^{cs}\times\{0\}$ which intersects $\cW^{u}_{loc}(f^n(x))$ in
a point $z$. The disk $D^n_x$ is the connected component
containing $f^{-n}(z)$ of $f^{-n}(D) \cap \Image B_x$.
\end{itemize}

Any family of such disks will converge to $\cD^{cs}_x$ by the
arguments sketched above (see the proof of Theorem 3.1 of
\cite{BF} for more details in a similar context).

To finish the proof is then enough to show that the plaques
$\cD^{cs}_x$ are coherent in the sense that if $y \in \cD^{cs}_x$
then $\cD^{cs}_x \cap \cD^{cs}_y$ is relatively open in
$\cD^{cs}_x$. To see that this holds in general we shall argue in
a similar way as in Lemma \ref{LemaRelacionEquivalencia} to take
advantage of the trapping property as well as the continuity of
the map $B$.

For each $x\in M$ we consider the set $\cW^{cs}_x$ defined as
$\bigcup_n f^{-n}(\cD^{cs}_{f^n(x)})$. Notice that $\cW^{cs}_x$ is
an immersed copy (in principle not injective) of $\RR^{\dim
E^{cs}}$ in $M$. We shall use in $\cW^{cs}_x$ the topology induced
by this immersion (i.e. the intrinsic topology and not the one
given as a subset of $M$).

To show that $E^{cs}$ is integrable, it is enough to show that
$\cW^{cs}_x$ is a partition of $M$ and that each leaf is
injectively immersed. Assume then that $y \in \cW^{cs}_x$, we must
show that $\cD^{cs}_y \en \cW^{cs}_x$. This will conclude since by
local uniqueness this gives that $\cW^{cs}_x$ is injectively
immersed and that the sets $\cW^{cs}_x$ are disjoint or coincide.

Consider, for $x \in M$ the sets $E_n = \{ y \in \cW^{cs}_x \ : \
f^n(\overline{\cD}^{cs}_y) \en \cD^{cs}_{f^n(x)} \}$. If one shows
that $\cW^{cs}_x = \bigcup_n E_n$ one completes the proof of the
Theorem since $f^n(\overline{\cD}^{cs}_y) \en \cD^{cs}_{f^n(x)}$
implies that $\cD^{cs}_y \en \cW^{cs}_x$ which as argued above
will imply that $\{\cW^{cs}_x \}_{x \in M}$ is an $f$-invariant
foliation tangent to $E^{cs}$.

The proof that the union $\bigcup_n E_n$ is closed in $\cW^{cs}_x$
is the same as in Lemma \ref{LemaRelacionEquivalencia}. The proof
of openness is slightly more delicate that in that case since we
do not know coherence in principle. However, coherence is easy to
establish for points which are nearby and that this is exactly
what we need to show to show openness.

To see that $E_n$ is open consider $z$ is close enough to $y$
verifying that $f^n(\overline{\cD}^{cs}_y) \en \cD^{cs}_{f^n(x)}$.
We must show that $z \in E_n$. To see this, it is enough to show
that given $y \in M$ one has that for $z$ in a small neighborhood
of $y$ in $\cD^{cs}_y$ it holds that $f(\cD^{cs}_z) \en
\cD^{cs}_{f(y)}$. By continuity of $B$ and the trapping property,
it follows that for $z$ in a neighborhood of $y$ in $\cD^{cs}_y$
the image by $f$ of $\Image B_z$ traverses the image of
$B_{f(y)}$. The characterization of $\cD^{cs}_z$ as limits of
disks as explained above implies that $f(\cD^{cs}_z) \en
\cD^{cs}_y$ as desired and concludes the proof.

%
%On the other hand, one can regard $\cD^{cs}_y$ as a points such
%that $f^k(y) \in \Image (B_{f^k(x)})$ for every $k\geq n$. Now, if
%$z \in \cW^{cs}_x$ is sufficiently close to $y$ one has that the
%set of points for which stay in the image of $B_{f^k(x)}$ for
%$k\geq n$ contains the $\cD^{cs}_z$ thanks to the trapping
%property and continuity of the function $B$. This concludes.

\lqqd

It  is natural to expect that for a partially hyperbolic
diffeomorphism $f: \TT^3 \to \TT^3$ with splitting $T\TT^3 =
E^{cs} \oplus E^u$ (with $\dim E^{cs}=2$) isotopic to a linear
Anosov automorphism with two-dimensional stable bundle, property
($\ast$) will be satisfied. To show this, one possibility would be
to show injectivity of the semiconjugacy to the linear model along
unstable manifolds but we have not succeed in doing so. A positive
answer would improve the results of \cite{Pot-PHFoliations} in
this context.

%%%%%%%%%%%%%%%%%%%%%%%%%%%%%%%%%%%%%%%%%%%%%%%%%%%%%%%%%%%%%%%%%%%%%%%%%%%%%
%%%%%%%%%%%%%%%%%%%%%%%%%%%%%%%%%%%%%%%%%%%%%%%%%%%%%%%%%%%%%%%%%%%%%%%%%%%%5
%%%%%%%%%%%%%%%%%%%%%%%%%%%%%%%%%%%%%%%%%%%%%%%%%%%%%%%%%%%%%%%%%%%%%%%%%55555

\appendix

\section{A non-trivial decomposition of the plane admitting
homotheties}\label{Appendix-Decomposition}

We shall denote as $d_2:\RR^2 \to \RR^2$ to the map

$$ d_2 (x) = \frac x 2 .$$

The goal of this appendix is to prove the following Theorem.

\begin{teo}\label{teoPPal} There exists a $C^\infty$-diffeomorphism $f:\RR^2 \to \RR^2$ and a constant $K>0$ such that the following properties are verified:
\bi \item[-] there exists a H\"{o}lder continuous cellular map
$h:\RR^2\to \RR^2$ such that $d_{C^0}(h,id) < K$ and $d_2 \circ h
= h \circ f$, \item[-] there exist open sets $V_1$ and $V_2$ such
that \bi \item $\overline{V_1}\cap \overline{V_2}=\emptyset$
    \item $h(V_i)= \RR^2$ for $i=1,2$.
    \item $f(\overline{V_i}) \en V_i$ for $i=1,2$.
\ei \item[-] the $C^{\infty}$ norm of $f$ and $f^{-1}$ is smaller
than $K$. \ei
\end{teo}

A direct consequence of this Theorem is the existence of $h:\RR^2
\to \RR^2$ whose fibers are all non trivial and cellular
(decreasing intersection of topological disks), the existence of
these decompositions of the plane had been shown by Roberts
\cite{RobertsDescomposiciones}.

\subsection{Construction of $f$.}\label{sectionConstruccion}

We start by considering a curve $\gamma= \{0\} \times [-\frac 1 4,
\frac 1 4]$.

Clearly, $\gamma \subset B_0 = B_1(0)$ the ball of radius one on
the origin. Consider also the sets $B_n = B_{2^n}(0)$ for every
$n\geq 0$. It follows that

$$ \RR^2= \bigcup_{n\geq 0} B_n $$

We shall define $f:\RR^2 \to \RR^2$ with the desired properties in
an inductive manner, first in $B_0$ and then in the annulus $B_n
\setminus B_{n-1}$ with arbitrary $n\geq 1$.
%
%Choose $K>0$ with the following properties:
%
%\bi \item[-] $K> 2diam (B_0) = 4$. \item[-] Every time there is a
%union $A_1 \cup A_2$ of two disjoint open sets in an annulus of
%the form $B_n \backslash B_{n-1}$ such that each of them is
%$K-$dense\footnote{We shall say a set $A$ is $a-$dense in $B$ with
%$a>0$ if for every point in $B$ there is a point in $A$ which is
%at distance strictly smaller than $a$ of it.} in $B_n \backslash
%B_{n-1}$ there exists a $C^\infty$-diffeomorphism $g:B_n
%\backslash B_{n-1} \to B_n \backslash B_{n-1}$ such that
%    \bi
%    \item $g$ coincides with the identity at the $K/10$ neighborhood of $\partial B_n$ and $\partial B_{n-1}$.
%    \item The $C^{0}$ distance between $g$ and the identity is smaller than $K$.
%    \item There exists a constant $a_K$ which depends only on $K$ and such that the $C^\infty$ norm of $g$ and its inverse is smaller than $a_K$.
%    \item The image by $g$ of $A_i$ is $K/2$-dense in $B_n \backslash B_{n-1}$ for every $i=1,2$.
%    \ei
%\ei
%
%To prove that such a $K$ exists, consider a function $\kappa:
%[diam (B_0), +\infty) \to \RR^+ \cup \{\infty\}$ such that
%$\kappa(t)$ is the infimum of the values which allow to find a
%$C^\infty$ function at $C^{\infty}$-distance smaller than
%$\kappa(t)$ as above which allows to make $t/2$-dense sets out of
%$t-$dense ones in the above context. It can be seen that $\kappa$
%is decreasing, so that it is enough to show that there exists $t$
%such that $\kappa(t) <\infty$ in order to show that $K$
%exists\marg{Corroborar}.

Define $f_0 : \overline{B_0} \to B_0$ a $C^\infty$ embedding and
open sets $V_1^0= \RR \times (-1/5,-1/6)$ and $V_2^0= \RR \times
(1/6,1/5)$ such that:

\bi \item[(a)] $f_0$ coincides with $d_2$ in a small neighborhood
of $\partial B_0$. \item[(b)] $\bigcap_{n\geq 0} f_0^n (B_0) =
\gamma$. \item[(c)] $f_0 (\overline{V_i^0})\en V_i^0$ for $i=1,2$.
\ei

Assume now that for some sufficiently large constants $K_1$ and
$K_2>0$ we have defined a $C^\infty$-diffeomorphism $f_n: B_n \to
B_{n-1}$ and disjoint open connected sets $V_1^n$ and $V_2^n$
(homeomorphic to a band $\RR \times (0,1)$) such that:

\bi \item[(I1)] $f_n|_{B_{n-1}}=f_{n-1}$ and $V_i^{n-1} \en V_i^n$
for $i=1,2$. \item[(I2)] the $C^\infty$-distance between $f_n$ and
$d_2$ in $B_n$ is smaller than $K_1$. \item[(I3)]
$f_n(\overline{V_i^n}) \en V_i^{n-1}$ for $i=1,2$ and
$f_n^n(\overline{V_i^n})$ disconnects $B_0$. \item[(I4)] $V_i^n$
contains balls of radius $\frac{1}{10}$ in every ball of radius
$\frac{K_2}{2}$ in $B_n$.
\item[(I5)] $f_n$ coincides with $d_2$ in
a $\frac{K_2}{10}$-neighborhood of $\partial B_n$. \ei

We must now construct $f_{n+1}$ and the sets $V_{i}^{n+1}$
assuming we had constructed $f_n$ and $V_i^n$.

To construct $f_{n+1}$ and $V_i^{n+1}$ we notice that in order to
verify (I1), it is enough to define $f_{n+1}$ in $B_n \backslash
B_{n-1}$ as well as to add to $V_i^n$ an open set in
$B_{n+1}\backslash B_{n}$ which verifies the desired hypothesis.

Consider $d_2^{-1}(V_i^n) \cap B_{n+1}\backslash B_n$. Since
$V_i^n$ satisfies property (I4) one has that $d_2^{-1}(V_i^n)$
contains a ball of radius $\frac{1}{5}$ in every ball of radius
$K_2$ of $B_{n+1}\backslash B_{n}$ for $i=1,2$.

Now, we consider a diffeomorphism $\varphi_n$ which is
$K_1-C^\infty$-close to the identity, coincides with the identity
in the $\frac{K_2}{10}$-neighborhoods of $\partial B_{n+1}$ and
$\partial B_n$ and such that $\varphi_n(V_i^n)$ contains a ball of
radius $\frac{1}{10}$ in every ball of radius $\frac{K_2}{2}$ of
$B_{n+1}$ for $i=1,2$. The existence of such $\varphi_n$ is
assured provided the value of $K_1$ is large enough with respect
to $K_2$.

We define then $f_{n+1}$ in $B_{n+1}\backslash B_n$ as $d_2 \circ
\varphi_n^{-1}$ which clearly glues together with $f_n$ and
satisfies properties (I2) and (I5).

To define $V_i^{n+1}$ we consider a very small $\eps>0$ (in order
that $\varphi_n(V_i^n)$ also verifies (I4)) and for each boundary
component $C$ of $\varphi_n(V_i^n)$ (which is a curve) we consider
a curve $C'$ which is at distance less than $\eps$ of $C$ inside
$\varphi_n(V_i^n)$ and such that each when it approaches $C\cap
\partial B_n$ the distance goes to zero and when it approaches $C
\cap\partial B_{n+1}$ the distance goes to $\eps$. This allows to
define new $V_i^{n+1}$ as the open set delimited by these curves
united with the initial $V_i^n$. It is not hard to see that it
will satisfy (I3) and (I4).

We have then constructed a $C^\infty$-diffeomorphism $f:\RR^2 \to
\RR^2$ which is at $C^\infty$ distance $K_1$ of $d_2$ and such
that there are two disjoint open connected sets $V_1$ and $V_2$
such that $f(\overline{V_i})\en V_i$. and such that both of them
are $\frac{K_2}{2}$-dense in $\RR^2$.

\subsection{Proof of the Theorem}

We first show the existence of a continuous function $h:\RR^2 \to
\RR^2$ conjugating $f$ to $d_2$ which is close to the identity.

This argument is quite classical: consider a point $x \in \RR^2$,
so, since $d_{C^0}(f,d_2)<K_1$ we get that the orbit $\{f^n(x)\}$
is in fact a $K_1-$pseudo-orbit for $d_2$. Since $d_2$ is
infinitely expansive, there exists only one orbit $\{d_2^n(y)\}$
which $\alpha(K_1)$-shadows $\{f^n(x)\}$ and we define $h(x)=y$
(in fact, in this case, it suffices with the past pseudo-orbit to
find the shadowing).

We get that $h$ is continuous since when $x_n \to x$ then the
pseudo-orbit which shadows must rest near for more and more time,
and then, again by expansivity, one concludes. This implies also
that $h$ is onto since it is at bounded distance of the identity.

Now, consider any ball $B$ of radius $100 \alpha(K)$ in $\RR^2$,
it is easy to see that $f(B)$ is contained in a ball of radius
$50\alpha(K)$ and then, we get a way to identify the preimage of
points by $h$. Consider a point $x \in \RR^2$, we get that

$$h^{-1}(h(x)) = \bigcap_{n>0} f^n( B_{100 \alpha(K)} (f^{-n}(x))) $$

\noindent which implies that $h$ is cellular.

It only remains to show that the image under $h$ of both $V_1$ and
$V_2$ is the whole plane.

\begin{lema} $h(V_i)= \RR^2$ for $i=1,2$.
\end{lema}

\dem We shall show that $h(\overline{V_i})$ is dense. Since it is
closed, this will imply that it is in fact the whole plane, and
using the semiconjugacy and the fact that $f(\overline{V_i})\en
V_i$ we get the desired property.

To prove that $h(\overline{V_i})$ is dense, we consider an
arbitrary open set $U\en \RR^2$. Now, choose $n_0$ such that
$d_2^{-n_0}(U)$ contains a ball of radius $10 \alpha(K)$. We get
that $h^{-1}(d_2^{-n_0}(U))$ contains a ball of radius $9
\alpha(K)$ and thus, since $\alpha(K)>K$, we know that since $V_i$
is $K/2$-dense, we get that $V_i \cap h^{-1}(d_2^{-n_0}(U)) \neq
\emptyset$. So, since $f(V_i) \en V_i$ we get that $V_i \cap
f^{n_0} \circ h^{-1}(d_2^{-n_0}(U)) \neq \emptyset$ which using
the semiconjugacy gives us that $h(V_i) \cap U \neq \emptyset$.

This concludes.

\lqqd

H\"{o}lder continuity of $h$ follows as in Theorem 19.2.1 of \cite{KH}
(see also \cite{Pot-Tesis}). Notice that the exponent of H\"{o}lder
continuity cannot be larger than $\frac{1}{2}$ since the boundary
of $V_i$ is sent as a space-filling curve.

\end{document}